\documentclass{article}
\usepackage{amsmath, amssymb}

\newtheorem{theorem}{Theorem}[section]
\newtheorem{proposition}[theorem]{Proposition}
\newtheorem{corollary}[theorem]{Corollary}

\newtheorem{lemma}[theorem]{Lemma}

\def\salt{\vspace{0.3 true cm}}

\def\n{\noindent}

\def\bop{\noindent {\bf Proof.}$ \;$ }
\def\eop{\hfill $\Box$ \vspace{0.5cm}}
\def\vareps{{\varepsilon}}
\def\conv{{\rm conv}}
\def\span{{\rm span}}

\def\cl{{\rm cl}}

\def\sgn{{\rm sgn}}

\def\bR{\mathbb{R}}

\def\cF{{\mathcal F}}

\def\cL{{\mathcal L}}
\def\cB{{\mathcal B}}

\def\cX{{\mathcal X}}

\def\cM{{\mathcal M}}
\def\cN{{\mathcal N}}

\def\cZ{{\mathcal Z}}
\def\bR{{\mathbb R}}

\def\bN{{\mathbb N}}

\begin{document}
\title{Bochner integrals and neural networks \thanks{appeared in {\it Handbook on Neural Information Processing}, Monica Bianchini, Marco Maggini, Lakhmi C. Jain, Eds., Springer, ISRL Vol. 49, 2013, Chap. 6, pp. 183--214}}
\author{Paul C. Kainen \and Andrew Vogt (1943--2021)}
\maketitle 

\begin{abstract}
A Bochner integral formula $f = \cB-\int_Y w(y) \Phi(y) \,d\mu(y)$ is derived
that presents a function $f$ in terms of weights $w$ and a parametrized
family of functions $\Phi(y)$, $y$ in $Y$.  
Comparison is made to pointwise formulations,
norm inequalities relating pointwise and Bochner integrals are established,
$G$-variation and tensor products are studied, and examples are presented.
\end{abstract}

\salt


\salt

\section{Introduction}

A neural network utilizes data to find a function consistent 
with the data and with further ``conceptual'' data such as desired smoothness,
boundedness, or integrability.  
The weights for a neural net and the functions embodied in the hidden units
can be thought of as determining a finite sum that approximates some
function. This finite sum is a kind of quadrature for an integral formula
that would represent the function exactly.
  
This chapter uses abstract analysis to investigate neural networks.
Our approach is one of {\it enrichment}: not only is summation replaced
by integration, but also numbers are
replaced by real-valued functions on an input set $\Omega$,
the functions lying in a function space $\cX$. The functions, in turn,
are replaced 
by $\cX$-valued measurable functions $\Phi$ on a measure space $Y$ of 
parameters.  
The goal is to understand approximation of functions by neural
networks so that one can make effective choices of the parameters
to produce a good approximation.

To achieve this, we utilize Bochner integration. The idea of applying this
tool to neural nets is in Girosi and Anzellotti \cite{gian93} and we developed
it further in Kainen and K\r{u}rkov\'{a} \cite{kavk09}.
Bochner integrals are now being used in the theory of support vector machines
and reproducing kernel Hilbert spaces; see the recent book by Steinwart 
and Christmann \cite{sc08}, which has an appendix of
more than 80 pages of material on operator theory and Banach-space-valued integrals.  Bochner integrals are also widely used in probability theory in
connection with stochastic processes of martingale-type;
see, e.g., \cite{ab82,cpmr07}.   The
corresponding functional
analytic theory may help to bridge the gap between probabilistic questions 
and deterministic ones, and may be well-suited for issues that
arise in approximation via neural nets.  

Training to replicate given numerical data 
does not give a useful neural network for the same
reason that parrots make poor conversationalists.  The phenomenon of
overfitting shows that achieving fidelity to data at all
costs is not desirable; 
see, e.g., the discussion on 
interpolation in our other chapter in this book (Kainen, K\accent23urkov\'{a}, 
and Sanguineti \cite{kks-h}).  In approximation, we try to find a function close to
the data that achieves desired criteria
such as sufficient smoothness, decay at infinity, etc.  Thus, a method
of integration which produces functions {\it in toto}
rather than numbers could be quite useful.

Enrichment has lately been utilized by
applied mathematicians to perform image analysis and even
to deduce global properties of sensor networks from local information.  For
instance, the Euler characteristic, ordinarily thought of as a discrete
invariant, can be made into a variable of
integration \cite{bg09}.  
In the case of sensor networks, such an analysis can lead to effective
computations in which theory determines a minimal set of sensors \cite{dsgh07}.

By modifying
the traditional neural net focus on training sets of data so that we get
to families of functions in a natural way,
we aim to achieve methodological insight. Such a framework may lead to
artificial neural networks capable of performing more sophisticated tasks.

The main result of this chapter is Theorem \ref{th:main}
which characterizes functions to be approximated in terms of
pointwise integrals and Bochner integrals, and provides inequalities
that relate corresponding norms.  The relationship between integral
formulas and neural networks has long been noted; e.g., 
\cite{it91,ba93,mhmi94,fg95-paris,ma96,vkkakr97}
We examine integral formulas in depth and extend their significance
to a broader context.

An earlier version of the Main Theorem, including the bounds
on variational norm by the $L^1$-norm of the weight function
in a corresponding integral formula, was given in \cite{kavk09}
and it also utilized {\it functional} (i.e., Bochner) integration.
However, the version here is more general and further shows that
if $\phi$ is a real-valued function on $\Omega \times Y$ (the cartesian
product of input and parameter spaces), then the associated map
$\Phi$ which maps the measure space to the Banach space defined by
$\Phi(y)(x) = \phi(x,y)$ is measurable; cf. \cite[Lemma 4.25, p. 125]{sc08}
where $\Phi$ is the ``feature map.''

Other proof techniques are available for parts of the Main Theorem.  
In particular, K\accent23urkov\'a \cite{vk09}
gave a different argument for part (iv) of the theorem, using a characterization of 
variation via peak functionals \cite{vksh98} as well as  
the theorem of Mazur (Theorem \ref{th:mazur})
used in the proof of Lemma \ref{pr:mvt}.
But the Bochner integral approach reveals some unexpected aspects of functional approximation which may be relevant for neural network applications. 

Furthermore, the treatment of analysis and topology utilizes a number of
basic theorems from the literature and provides an introduction
to functional analysis motivated by its applicability.  This is a case where
neural nets provide a fresh perspective on classical mathematics.  
Indeed, theoretical results proved here were obtained
in an attempt to better understand neural networks.

An outline of the paper is as follows: In section 2 we discuss variational
norms; sections 3 and 4 present needed material on Bochner integrals.  The
Main Theorem (Theorem \ref{th:main}) on integral formulas is given in Section
5.  In section 6 we show how to apply the Main
Theorem to an integral formula for the Bessel
potential function in terms of Gaussians.  In section 7 we show how this leads
to an inequality involving Gamma functions and provide an alternative
proof by classical means.  Section 8 interprets and extends the Main Theorem
in the language of tensor products.  Using tensor products, we replace
individual $\cX$-valued
$\Phi$'s by families $\{\Phi_j: j \in J\}$ of such functions.
This allows more nuanced representation of the function to be approximated.
In section 9 we give a detailed example
of concepts related to $G$-variation, while section 10 considers the 
relationship between pointwise integrals and evaluation
of the corresponding Bochner integrals.  Remarks on future
directions are in section 11, and the chapter concludes
with two appendices and references.

\section{Variational norms and completeness}

We assume that the reader has a reasonable acquaintance with functional
analysis but have attempted to keep this chapter self-contained.
Notations and basic definitions are given in Appendix I, while 
Appendix II has the precise statement of several important theorems
from the literature which will be needed in our development.

Throughout this chapter, all linear spaces are over the reals $\bR$.
For $A$ any subset of a linear
space $X$, $b \in X$, and $r \in \bR$, 
$$b + rA := \{b + ra \,|\, a \in A \} = \{y \in X : y = b + ra, a \in A\}.$$
Also, we sometimes use the abbreviated notation  

\begin{equation}
\|\cdot\|_1 = \|\cdot\|_{L^1(Y,\mu)} \;\;\mbox{and} 
\;\; \|\cdot\|_\infty = \|\cdot\|_{L^\infty(Y,\mu;\cX)};
\label{eq:notation}
\end{equation}

\n the standard notations on the right are explained in sections 12 and 4,
resp.  The symbol ``$\ni$'' stands for ``such that.''

A set $G$ in a normed linear space $\cX$ is {\it fundamental} (with
respect to $\cX$) if $\cl_{\cX} \; (\span \; G) = \cX$,
where closure depends only on the topology induced by the norm.
We call $G$ {\it bounded} with respect to $\cX$ if
$$s_{G,\cX} := \sup_{g \in G} \|g\|_{\cX} < \infty.$$

We now review $G$-{\it variation} norms.   
These norms, which arise in connection with approximation of functions,
were first considered by Barron \cite{ba92}, \cite{ba93}. He treated a
case where $G$ is a family of characteristic functions of sets satisfying a 
special condition.  The general concept, formulated by
K\r{u}rkov\'{a} \cite{vk97}, has been developed in such papers as
\cite{vksa02,vk03,vk08,JOTA09,JOTA10a,gs11}.

Consider the set

\begin{equation}
B_{G,\cX} := \cl_{\cX} \; (\conv( \; \pm G)), \mbox{where} \; \pm G
:= G \cup -G. \label{eq:BGX}
\end{equation}

\n
This is a symmetric, closed, convex subset of $\cX$, with
Minkowski functional
$$\|f\|_{G,\cX} := \inf \{\lambda > 0: f/\lambda \in B_{G,\cX} \}.$$
The subset $\cX_G$ of $\cX$ on which this functional 
is finite is given by
$$\cX_G := \{f \in \cX: \exists \lambda > 0 \, \ni \,
f/\lambda \in B_{G,\cX} \}.$$ 
If $G$ is bounded, then $\|\cdot\|_{G,\cX}$ is a norm on
$\cX_G$.  In general $\cX_G$ may be a proper subset of $\cX$ even if $G$
is bounded and fundamental w.r.t. $\cX$. See the example at the end of
this section.  The inclusion
$\iota: \cX_G \subseteq \cX$ is linear and for
every $f \in \cX_G$

\begin{equation}
\|f\|_{\cX} \leq \|f\|_{G,\cX} \,s_{G,\cX} \label{eq:basicIneq}
\end{equation}

\n Indeed, if $f/\lambda \in B_{G,\cX}$, then
$f/\lambda$ is a convex combination of elements of $\cX$-norm
at most $s_{G,\cX}$, so
$\|f\|_{\cX} \leq \lambda \, s_{G,\cX}$
establishing (\ref{eq:basicIneq})
by definition of variational norm.
Hence, if $G$ is bounded in $\cX$, the operator $\iota$ is
bounded with operator norm not exceeding $s_{G,\cX}$.

\begin{proposition}

Let nonempty $G \subseteq \cX$ a normed linear space.  Then\\
\\
(i) $\span \,G \subseteq \cX_G \subseteq \cl_\cX\, \span \,G$;\\
\\
(ii) $G$ is fundamental if and only if $\cX_G$ is dense in $\cX$;\\
\\
(iii) For $G$ bounded and $\cX$ complete,  
$(\cX_G, \|\cdot\|_{G,\cX})$ is a Banach space.
\label{pr:varban}
\end{proposition}

\bop (i) Let $f \in \span \,G$, then $f = \sum_{i=1}^n a_i g_i$,
for real numbers $a_i$ and $g_i \in G$.  We assume the $a_i$ are
not all zero since $0$ is in $\cX_G$.  Then 
$f = \lambda \sum_{i=1}^n |a_i|/\lambda (\pm g_i)$, where
$\lambda = \sum_{i=1}^n |a_i|$. Thus, $f$ is in 
$\lambda \conv (\pm G) \subseteq \lambda B_{G,\cX}$.
So $\|f\|_{G,\cX} \leq \lambda$ and $f$ is in $\cX_G$.

Likewise if $f$ is in $\cX_G$, then for some $\lambda > 0$, 
$f/\lambda$ is in 
$$B_{G,\cX} = \cl_\cX(\conv(\pm G)) \subseteq \cl_X(\span(G)),$$
so $f$ is in $\cl_X(\span(G))$.

(ii) Suppose $G$ is fundamental.  Then  
$\cX = cl_\cX(\span \,G) = \cl_\cX(\cX_G)$ by part (i). Conversely,
if $\cX_G$ is dense in $\cX$, then
$\cX = \cl_\cX(\cX_G) \subseteq \cl_\cX(\span \,G) \subseteq \cX$,
and $G$ is fundamental.

(iii)
Let $\{f_n\}$ be a Cauchy sequence in $\cX_G$. By
(\ref{eq:basicIneq})  $\{f_n\}$ is a Cauchy sequence in $\cX$ and
has a limit $f$ in $\cX$. The sequence $\|f_n\|_{G,\cX}$ is bounded
in $\cX_G$, that is, there is a positive number M such that for all
n $f_n/M \in B_{G,\cX}$. Since $B_{G,\cX}$ is closed in $\cX$, $f/M$
is also in $B_{G,\cX}$. Hence $\|f\|_{G,\cX} \leq M$ and $f$ is in
$\cX_G$. Now given $\epsilon > 0$ choose a positive integer $N$ such
that $\|f_n - f_k\|_{G,\cX} < \epsilon$ for $n, k \geq N$. In
particular fix $n \geq N$, and consider a variable integer $k \geq
N$. Then $\|f_k - f_n\|_{G,\cX} < \epsilon$. So $(f_k -
f_n)/\epsilon \in B_{G,\cX}$, and $f_k \in f_n + \epsilon B_{G,\cX}$
for all $k \geq N$. But $f_n + \epsilon B_{G,\cX}$ is closed in
$\cX$. Hence $f \in f_n + \epsilon B_{G,\cX}$, and $\|f -
f_n\|_{G,\cX} \leq \epsilon$. So the sequence converges to
$f$ in $\cX_G$.
\eop

 The following example illustrates several of the above concepts. 
 Take $\cX$ to be a real
 separable Hilbert space with orthonormal basis $\{e_n : n = 0, 1,
 ...\}$. Let $G = \{e_n : n = 0, 1, ...\}$. Then
 $$B_{G,\cX} = 
 \left \{ \sum_{n \geq 1}
 c_ne_n - \sum_{n \geq 1} d_ne_n\; : \forall n,\; c_n \geq 0, d_n \geq 0,
   \sum_{n \geq 1}(c_n + d_n) = 1 \right \}.$$

\n Now $f \in \cX$ is of the form $\sum_{n \geq 1} a_ne_n$
where $\|f\|_{\cX} = \sqrt{\sum_{n \geq 1} a_n^2}$, and if $f \in
\cX_G$, then $a_n = \lambda (c_n - d_n)$ for all $n$ and suitable
$c_n, d_n$.  The minimal $\lambda$ can
be obtained by taking $a_n = \lambda c_n$ when $a_n \geq 0$, and
$a_n = -\lambda d_n$ when $a_n < 0$. It then follows that
$\|f\|_{G,\cX} = \sum_{n \geq 1} |a_n|$. Hence when $\cX$ is
isomorphic to $\ell_2$, $\cX_G$ is isomorphic to
$\ell_1$.  As $G$ is fundamental, by part(ii) above, the
closure of  $\ell_1$ in $\ell_2$ is $\ell_2$.
This provides an example where $\cX_G$ is not
a closed subspace of $\cX$ and so, while it is a Banach
space w.r.t. the variational norm, it is not complete in
the ambient-space norm.\\

\section{Bochner integrals}

The Bochner integral replaces numbers with functions and represents
a broadranging extension,
generalizing the Lebesgue integral from real-valued
functions to functions with values in an arbitrary Banach space.  
Key definitions and theorems are summarized here for convenience,
following the treatment in \cite{za61} (cf. \cite{lc85}).  Bochner
integrals are used here  (as in \cite{kavk09}) in order to prove a
bound on variational norm.

Let $(Y, \mu)$ be a measure space.  Let $\cX$ be a Banach space
with norm $\|\cdot\|_{\cX}$. A function $s: Y \to \cX$ is  {\it simple} if it
has a finite set of nonzero values $f_j \in \cX$, 
each on a measurable subset $P_j$ of $Y$ with $\mu(P_j) < \infty$,
$1 \leq j \leq m$, and the $P_j$ are pairwise-disjoint.
Equivalently, a function $s$ is simple if it can be written in the 
following form:
\begin{equation} 
\hspace{1 in} s = \sum_{j=1}^m \kappa(f_j) \chi_{P_j} ,
\label{eq:sf}\end{equation}
where $\kappa(f_j):Y \to \cX$ denotes the
constant function with value $f_j$ and $\chi_{P}$ denotes the
characteristic function of a subset $P$ of $Y$.  
This decomposition is nonunique and we identify two
functions if they agree $\mu$-almost everywhere - i.e., the 
subset of $Y$ on which they disagree has $\mu$-measure zero.  

Define an $\cX$-valued function $I$ on the simple functions 
by setting for $s$ of form (\ref{eq:sf})
$$I(s,\mu) := \sum_{j=1}^m  \mu(P_j) f_j \in \cX.$$
This is independent of the decomposition of $s$ \cite[pp.130--132]{za61}.
A function $h:Y \to \cX$ is  {\it strongly measurable} (w.r.t. $\mu$) if there exists a sequence $\{ s_k\}$  of
simple functions such that for $\mu$-a.e. $y \in Y$
$$\lim_{k \to \infty} \|s_k(y) - h(y)\|_{\cX} = 0.$$

A function
$h:Y \to \cX$ is {\it Bochner integrable} (with respect to $\mu$)
if it is strongly measurable and there exists a sequence $\{s_k\}$
of simple functions $s_k: Y \to \cX$ such that
\begin{equation}
\hspace{1 in} \lim_{k \to \infty} \int_Y \|s_k(y) - h(y)\|_{\cX} d\mu(y) = 0.
\label{eq:bis} \end{equation}
If $h$ is strongly measurable and (\ref{eq:bis}) holds, then 
the sequence $\{I(s_k,\mu)\}$ is Cauchy and by completeness converges to an
element  in $\cX$.  This element, which is independent of the 
sequence of simple functions satisfying (\ref{eq:bis}), is called
the {\it Bochner integral of} $h$ (w.r.t. $\mu$) and denoted
$$I(h,\mu) \;\;\; \mbox{or} \;\;\; \cB-\int_Y h(y) d\mu(y).$$

Let $\cL^1(Y,\mu;\cX)$ denote the linear space of all 
strongly measurable functions from $Y$
to $\cX$ which are  Bochner integrable w.r.t. $\mu$;
let $L^1(Y,\mu;\cX)$ be the corresponding set of equivalence
classes (modulo $\mu$-a.e. equality).  It is easily shown that
equivalent functions have the same Bochner integral.
Then the following elegant characterization holds.

\begin{theorem}[Bochner]
Let $(\cX,\|\cdot\|_{\cX})$ be a Banach space and $(Y, \mu)$ a
measure space. Let $h:Y \to \cX$ be strongly measurable. Then
$$h \in \cL^1(Y,\mu;\cX) \mbox{ if and only
if } \int_Y \|h(y)\|_{\cX} d\mu(y) < \infty.$$\label{th:bo}
\end{theorem}

\n A consequence of this theorem is that 
$I: L^1(Y,\mu;\cX) \to \cX$ is a continuous linear operator
and
\begin{equation}
\hspace{-.3 in}\|I(h,\mu)\|_\cX = \left \| \cB-\int_Y h(y) \,d\mu(y) \right \|_\cX
\leq \|h\|_{L^1(Y,\mu;\cX)} := \int_Y \|h(y)\|_X d\mu(y).
\label{ineq:boch}
\end{equation}
In particular, the Bochner norm of $s$, $\|s\|_{L^1(Y,\mu;\cX)}$,
is $\sum_i \mu(P_i) \|g_i\|_\cX$, where $s$ is a simple function
satisfying (\ref{eq:sf}).

For $Y$ a measure space and $\cX$ a Banach space,
$h:Y \to \cX$ is {\it weakly measurable} if for every
continuous linear 
functional $F$ on $X$ the composite real-valued function
$F \circ h$ is measurable  \cite[pp. 130--134]{yo65}.
If $h$ is measurable, then it is weakly measurable
since measurable followed by continuous is measurable:
for $U$ open in $\bR$, 
$(F \circ h)^{-1}(U) = h^{-1}(F^{-1}(U))$.

Recall that a topological space
is {\it separable} if it has a countable dense subset.
Let $\lambda$ denote Lebesgue measure on $\bR^d$ and let $\Omega \subseteq \bR^d$
be $\lambda$-measurable,  $d \geq 1$.
Then $L^q(\Omega,\lambda)$ is separable when $1 \leq q < \infty$; e.g., 
\cite[pp. 208]{mb59}.
A function $h:Y \to \cX$ is {\it $\mu$-almost separably valued} ($\mu$-a.s.v.) if there exists a $\mu$-measurable subset
$Y_0 \subset Y$ with $\mu(Y_0)=0$ and $h(Y \setminus Y_0)$ is
a separable subset of $\cX$.

\begin{theorem}[Pettis]
Let $(\cX,\|\cdot\|_{\cX})$ be a Banach space and $(Y, \mu)$ a
measure space. Suppose $h:Y \to \cX$.  Then $h$ is strongly measurable
if and only if $h$ is weakly measurable and $\mu$-a.s.v.
\label{th:pe}
\end{theorem}

The following basic result (see, e.g., \cite{sb33}) was later extended 
by Hille to the more general class of closed operators.  But we only
need the result for bounded linear functionals, in which case the 
Bochner integral coincides with ordinary integration.

\begin{theorem}
Let $(Y,\nu)$ be a measure space, let $\cX$, $\cX'$ be Banach
spaces, and let $h \in \cL^1(Y,\nu;\cX)$. If $T: \cX \to \cX'$ is a 
bounded linear operator, then $T \circ h \in \cL^1(Y,\nu;\cX')$ and
$$T\left( \cB-\int_Y h(y) \,d\nu(y) \right) = \cB-\int_Y (T \circ h)(y) \,d\nu(y).$$
\label{pr:blf}\end{theorem}

There is a mean-value 
theorem for Bochner integrals (Diestel and Uhl
\cite[Lemma 8, p. 48]{du77}).  We give their
argument with a slightly clarified reference to the Hahn-Banach theorem.

\begin{lemma}
Let $(Y,\nu)$ be a finite measure space, let $X$ be a Banach space, and
let $h:Y \to \cX$ be Bochner integrable w.r.t. $\nu$.  
Then $$\cB-\int_Y h(y) \,d\nu(y) \;\in\; \nu(Y) \;\cl_X(\conv(\{\pm h(y): y \in Y\}).$$
\label{pr:mvt}
\end{lemma}

\bop
Without loss of generality, $\nu(Y) = 1$.  
Suppose  $f := I(h,\nu) \notin \cl_X(\conv(\{\pm h(y): y \in Y\})$.
By a consequence of the 
Hahn-Banach theorem 
given as Theorem \ref{th:mazur} in Appendix II below),
there is a continuous linear functional $F$ on $X$ such that $F(f) >  
\sup_{y \in Y} F(h(y))$. Hence, by Theorem \ref{pr:blf},
$$\sup_{y \in Y} F(h(y)) \geq 
\int_Y F(h(y)) d\nu(y) = F(f) > \sup_{y \in Y} F(h(y)).$$
which is absurd.  \eop

\section{Spaces of Bochner integrable functions}

In this section, we derive a few consequences of the results from
the previous section which we shall need below.


A measurable function $h$ from
a measure space $(Y,\nu)$ to a normed linear space $\cX$ is called 
{\it essentially bounded} (w.r.t. $\nu$) if there exists a $\nu$-null set 
$N$ for which $$\sup_{y \in Y \setminus N} \|h(y)\|_\cX < \infty.$$ Let 
$\cL^{\infty}(Y,\nu; \cX)$ denote the linear space of all 
strongly measurable, essentially
bounded functions from $(Y,\nu)$ to $\cX$. Let $L^\infty(Y, \nu; \cX)$ be
its quotient space mod the relation of equality $\nu$-a.e. 
This is a Banach space with norm
$$\|h\|_{L^{\infty}(Y,\nu;X)} := \inf \{B \geq 0: \exists \; \nu \mbox{-null } 
N \subset Y \; \ni \;\|h(y)\|_X \leq B, \; \forall y \in Y \setminus N \}.$$
To simplify notation, we sometimes write $\|h\|_\infty$ for
$\|h\|_{L^{\infty}(Y,\nu;X)}$
Note that if $\|h\|_\infty = c$, then $\|h(y)\|_\cX \leq c$
for $\nu$-a.e. $y$.  Indeed, for positive integers $k$,
$\|h(y)\|_\cX \leq c + (1/k)$ for $y$ not in a set of measure zero $N_k$
so $\|h(y\|_\cX \leq c$ for $y$ not in the union $\bigcup_{k\geq 1} N_k$
also a set of measure zero.

We also have a useful fact whose proof is immediate.

\begin{lemma}
For every measure space $(Y,\mu)$ and Banach space $\cX$, 
the natural map $\kappa_\cX: \cX \to L^\infty(Y,\mu;\cX)$ 
associating to each element $g \in \cX$ the constant function 
from $Y$ to $\cX$ given by $(\kappa_\cX(g))(y) \equiv g$ for all $y$ in $Y$ 
is an isometric linear embedding.
\end{lemma}

\begin{lemma}
Let $\cX$ be a separable Banach space, let $(Y,\mu)$ be a measure space,
and let $w:Y \to \bR$ and $\Psi:Y \to \cX$ be $\mu$-measurable functions.
Then $w \Psi$ is strongly measurable.
\label{lm:sm}
\end{lemma}

\bop
By definition, $w \Psi$ is the function from $Y$ to $\cX$ defined by
$$w \Psi: \; y \mapsto w(y) \Psi(y),$$
where the multiplication is that of a Banach space element by a real number.
Then $w \Psi$ is measurable because it is obtained from a pair of
measurable functions by applying scalar multiplication which is continuous.
Hence, by separability, Pettis' Theorem \ref{th:pe}, and the
fact that measurable implies weakly measurable, we have
strong measurability for $w \Psi$ (cf. \cite[Lemma 10.3]{lc85}).
\eop

If $(Y,\nu)$ is a finite measure space, $\cX$ is a Banach space, 
and $h:Y \to \cX$ is 
strongly measurable and essentially bounded, then $h$ is Bochner
integrable by Theorem \ref{th:bo}.  The following lemma, which follows from 
Lemma \ref{lm:sm}, allows us to 
weaken the hypothesis on the function by further constraining the space $\cX$.

\begin{lemma}
Let $(Y,\nu)$ be a finite measure space, $\cX$ a separable Banach space,
and $h:Y \to \cX$ be $\nu$-measurable and essentially bounded w.r.t. $\nu$.
Then $h \in \cL^1(Y,\nu; \cX)$ and 
$$\int_Y \|h(y)\|_\cX d\nu(y) \leq \nu(Y) \|h\|_{L^\infty(Y,\nu; \cX)}.$$
\label{lm:bi}
\end{lemma}

Let $w \in \cL^1(Y,\mu)$, and
let $\mu_w$ be defined for $\mu$-measurable 
$S \subseteq Y$ by $\mu_w(S) := \int_S |w(y)| d\mu(y)$.
For $t \neq 0$, $\sgn(t) := t/|t|$.

\begin{theorem}
Let $(Y,\mu)$ be a measure space, 
$\cX$ a separable Banach space; let $w \in \cL^1(Y,\mu)$ be nonzero
$\mu$-a.e., let $\mu_w$ be the measure defined above, and let
$\Phi:Y \to \cX$ be $\mu$-measurable.  If one of
the Bochner integrals 
$$\cB-\int_Y w(y) \Phi(y) d\mu(y), \; \;\;  
\cB-\int_Y \sgn(w(y)) \Phi(y) d\mu_w(y)$$
exists, then both exist and are equal.
\label{th:boeq}
\end{theorem}

\bop
By Lemma \ref{lm:sm}, both $w \Phi$ and $(\sgn \circ w) \Phi$ are strongly measurable.
Hence, by Theorem \ref{th:bo}, the respective Bochner integrals 
exist if and only
if the $\cX$-norms of the respective integrands have finite 
ordinary integral. But 
\begin{equation}
\int_Y \| [(\sgn \circ w)\Phi ](y)\|_{\cX} d\mu_w(y) 
= \int_Y \|w(y) \Phi(y)\|_{\cX} d\mu(y),
\label{eq:same}
\end{equation}
so the Bochner integral $I((\sgn \circ w)\Phi, \mu_w)$ exists exactly 
when $I(w \Phi, \mu)$ does.
Further, the respective Bochner integrals are equal since for 
any continuous linear functional $F$ in $\cX^*$, by Theorem \ref{pr:blf}
\begin{eqnarray*}
F\left(\cB-\int_Y w(y) \Phi(y) d\mu(y) \right) 
= \int_Y F(w(y) \Phi(y)) d\mu(y) \\
= \int_Y w(y) F(\Phi(y)) d\mu(y)
= \int_Y sgn(w(y))|w(y)| F(\Phi(y)) d\mu(y)\\
= \int_Y sgn(w(y)) F(\Phi(y)) d\mu_w(y)
= \int_Y F(sgn(w(y)) \Phi(y)) d\mu_w(y) \\
= F\left(\cB-\int_Y sgn(w(y)) \Phi(y) d\mu_w(y)\right).
\label{eq:fin-boch}
\end{eqnarray*}
\eop

\begin{corollary}
Let $(Y,\mu)$ be a $\sigma$-finite measure space, 
$\cX$ a separable Banach space, $w: Y \to \bR$ be in $\cL^1(Y,\mu)$
and $\Phi:Y \to \cX$ be in $\cL^\infty(Y,\mu; \cX)$.  Then $w \Phi$ is
Bochner integrable w.r.t. $\mu$.
\label{co:bi}
\end{corollary}

\bop
By Lemma \ref{lm:sm}, $(\sgn \circ w) \Phi$ is strongly measurable,
and Lemma \ref{lm:bi} then implies that the Bochner integral
$I((\sgn \circ w) \Phi, \mu_w)$ exists since $\mu_w(Y) = \|w\|_{L^1(Y,\mu)} < \infty$.
So $w \Phi$ is Bochner integrable by Theorem \ref{th:boeq}.
\eop

\section{Main theorem}

In the next result, we show that certain types of integrands yield
integral formulas for functions $f$ in a Banach space of $L^p$-type both
pointwise and at the level of Bochner integrals.  Furthermore, the variational norm of $f$ is shown to be bounded
by the $L^1$-norm of the weight function from the integral formula.
Equations (\ref{eq:th2}) and (\ref{eq:fbi}) and part (iv) of this theorem were
derived in a similar fashion by one of us with K\r{u}rkov\'{a} in \cite{kavk09}
under more stringent hypotheses; see also \cite[eq. (12)]{fg95-paris}.

\begin{theorem}
Let $(\Omega,\rho)$, $(Y,\mu)$ be $\sigma$-finite measure spaces, 
let $w$ be in $\cL^1(Y,\mu)$, let  
$\cX = L^q(\Omega,\rho)$, $q \in [1,\infty)$, be separable, let 
$\phi: \Omega \times Y \to \bR$ be $\rho \times \mu$-measurable, let
$\Phi: Y \to \cX$ be defined for each $y$ in $Y$ by 
$\Phi(y)(x) := \phi(x,y)$ for $\rho$-a.e. $x \in \Omega$
and suppose that for some $M < \infty$, 
$\|\Phi(y)\|_\cX \leq M$ for $\mu$-a.e. $y$.
Then the following hold:\\

(i) For $\rho$-a.e. $x \in \Omega$, the integral
$\int_Y w(y) \phi(x,y)d\mu(y)$ exists and is finite.\\

(ii) The function $f$ defined by
\begin{equation}
\;\;\; f(x) = \int_Y w(y) \phi(x,y)d\mu(y)
\label{eq:th} 
\end{equation}
is in $\cL^q(\Omega, \rho)$ and its equivalence class, also denoted
by $f$, is in $L^q(\Omega, \rho) = \cX$ and satisfies
\begin{equation}
\;\;\; \|f\|_\cX \leq \|w\|_{L^1(Y,\mu)}\;M.
\label{eq:th2} 
\end{equation}

(iii) The function $\Phi$ is measurable and hence in $\cL^\infty(Y,\mu;\cX)$,
and $f$ is the Bochner integral of $w \Phi$ w.r.t. $\mu$, i.e.,
\begin{equation}
f = \cB-\int_Y (w \Phi)(y) d\mu(y).
\label{eq:fbi}
\end{equation}

(iv) For $G = \{\Phi(y) : \|\Phi(y)\|_\cX \leq 
\|\Phi\|_{L^\infty(Y,\mu; \cX)} \}$, $f$ is in $\cX_G$, and
\begin{equation}
\|f\|_{G,\cX} \leq \|w\|_{L^1(Y,\mu)}
\label{ineq:varL1}
\end{equation}
and as in (\ref{eq:notation})
\begin{equation}
\|f\|_\cX \leq
\|f\|_{G,\cX} s_{G,\cX} \leq \|w\|_1 \|\Phi\|_\infty.\;\;\;
\label{ineq:thm}
\end{equation}
\label{th:main}
\end{theorem}

\bop
(i) Consider the function $(x, y) \longmapsto |w(y)||\phi(x,y)|^q$. This
is a well-defined $\rho \times \mu$-measurable function on $\Omega
\times Y$. Furthermore its repeated integral

$$\int_Y \int_{\Omega} \, |w(y)||\phi(x,y)|^q  d\rho(x)  d\mu(y)$$

\n exists and is bounded by $\|w\|_1 M^q$ since $\Phi(y) \in
L^q(\Omega,\rho)$ and $\|\Phi(y)\|_q^q \leq M^q$ for a. e. y. and $w
\in  L^1(Y,\mu)$. By Fubini's Theorem \ref{th:fubini} the
function $y \longmapsto |w(y)||\phi(x,y)|^q$ is in $L^1(Y,\mu)$ for
a.e. x. But the inequality

$$|w(y)||\phi(x,y)| \leq \max\{|w(y)||\phi(x,y)|^q,|w(y)|\} \leq
(|w(y)||\phi(x,y)|^q + |w(y)|) $$

\n shows that the function $y \longmapsto |w(y)||\phi(x,y)|$ is
dominated by the sum of two integrable functions. Hence the
integrand in the definition of $f(x)$ is integrable for a. e. x, and
$f$ is well-defined almost everywhere.

(ii) The function $G(u) = u^q$ is a convex function for $u \geq 0$.
Accordingly by Jensen's inequality (Theorem \ref{th:jensen} below),

$$ G \left (\int_Y |\phi(x,y)| d\sigma(y) \right ) \leq \ \int_Y
G(|\phi(x,y)|)d\sigma(y)$$

\n provided both integrals exist and $\sigma$ is a probability
measure on the measurable space $Y$.  We take $\sigma$ to be defined
by the familiar formula:

$$ \sigma(A) = \frac{\int_A |w(y)| d\mu(y)}{\int_Y |w(y)|
d\mu(y)}$$

\n for $\mu$-measurable sets A in Y, so that integration with
respect to $\sigma$ reduces to a scale factor times integration of
$|w(y|d\mu(y)$.
Since we have established that both $|w(y)||\phi(x,y)|$ and
$|w(y)||\phi(x,y)|^q$ are integrable with respect to $\mu$ for a.e.
x, we obtain:

$$ |f(x)|^q \leq \|w\|_1^q G(\int_Y |\phi(x,y)| d\sigma(y)) \leq
\|w\|_1^q \int_Y G(|\phi(x,y)|)d\sigma(y)$$
$$\;\;\; = \|w\|_1^{q-1}\int_Y
|w(y)||\phi(x,y)|^q d\mu(y) $$

\n for a.e. x. But we can now integrate both side with respect to
$d\rho(x) $ over $\Omega$ because of the integrability noted above in
connection with Fubini's Theorem. Thus $f \in \cX =
L^q(\Omega,\rho)$ and $\|f\|^q_\cX \leq \|w\|_1^q \, M^q$, again
interchanging order.

(iii) First we show that $\Phi^{-1}$ of the open ball centered at $g$ of radius $\vareps$,  
$B(g,\vareps) := \{y: \|\Phi(y) - g\|_\cX < \vareps \}$, is a 
$\mu$-measurable subset
of $Y$ for each $g$ in $\cX$ and $\vareps > 0$.  Note that
$$\|\Phi(y) - g\|_\cX^q = \int_{\Omega} |\phi(x,y) - g(x)|^q d\rho(x)$$
\n for all $y$ in $Y$ where $x \mapsto \phi(x,y)$ and $x \mapsto g(x)$ 
are $\rho$-measurable functions representing the elements $\Phi(y)$ and $g$
belonging to $\cX = L^q(Y,\mu)$.  Since $(Y,\mu)$ is $\sigma$-finite, we
can find a strictly positive function $w_0$ in $\cL^1(Y,\mu)$. (For example,
let $w_0 = \sum_{n \geq 1} (1/n^2) \chi_{Y_n}$, where $\{Y_n: n \geq 1 \}$
is a countable disjoint partition of $Y$ into $\mu$-measurable sets of finite
measure.)  Then $w_0(y) |\phi(x,y) - g(x)|^q$ is a $\rho \times \mu$-measurable
function on $\Omega \times Y$, and 
$$\int_Y \int_{\Omega} w_0(y) |\phi(x,y) - g(x)|^q d\rho(x) d\mu(y) \leq
\|w_0\|_{L^1(Y,\mu)} \vareps^q.$$
By Fubini's Theorem \ref{th:fubini}, 
$y \mapsto w_0(y) \|\Phi(y) - g\|_\cX^q$ is $\mu$-measurable.  Since $w_0$ is $\mu$-measurable and strictly positive,
$y \mapsto \|\Phi(y) - g\|_\cX^q$ is also $\mu$-measurable and so 
$B(g,\vareps)$ is measurable.  Hence, $\Phi: Y \to \cX$ is measurable.
Thus, $\Phi$ is essentially bounded, with essential sup
$\|\Phi\|_{L^\infty(Y,\mu;\cX)} \leq M$.  (In (\ref{eq:th2}), $M$ can be
replaced by this essential sup.)

By Corollary \ref{co:bi}, $w \Phi$ is Bochner integrable.
To prove that $f$ is the Bochner integral, using Theorem \ref{th:boeq},
we show that for each bounded linear functional 
$F \in \cX^*$, $F(I(\sgn \circ w \Phi,\mu_w)) = F(f)$. By
the Riesz representation theorem \cite[p. 316]{masa01}, for any such $F$
there exists a (unique) $g_F \in \cL^{p}(\Omega,\rho)$, $p = 1/(1-q^{-1})$, 
such that for
all $g \in \cL^q(\Omega,\rho)$, $F(g) = \int_{\Omega} g_F(x) g(x)
d\rho(x)$. 
By Theorem \ref{pr:blf}, 
$$F(I((\sgn \circ w) \Phi,\mu_w)) = \int_Y F\left(\sgn(w(y)) \Phi(y)\right) d\mu_w(y).$$
But for $y \in Y$, $F(\sgn(w(y)) \Phi(y)) = \sgn(w(y)) F(\Phi(y)$, so
$$F(I((\sgn \circ w) \Phi,\mu_w)) = 
\int_Y \int_{\Omega} w(y) g_F(x) \phi(x,y) d\rho(x) d\mu(y).$$
Also, using (\ref{eq:th}),
$$F(f) = \int_{\Omega} g_F(x) f(x) d\rho(x) = \int_{\Omega} \int_Y w(y)
g_F(x) \phi(x,y) d\mu(y) d\rho(x).$$ 
The integrand of the iterated integrals is measurable with 
respect to the product measure $\rho \times \mu$, so by 
Fubini's Theorem   the iterated integrals are equal 
provided that one of the corresponding absolute integrals is finite.  Indeed,
\begin{equation}
\int_Y \int_{\Omega} |w(y) g_F(x) \phi(x,y)| d\rho(x) d\mu(y) 
= \int_Y \|g_F \Phi(y) \|_{L^1(\Omega,\rho)} d\mu_w(y).\;\;
\label{eq:fub}
\end{equation}
By H\"{o}lder's inequality, for every $y$, 
$$\|g_F \Phi(y) \|_{L^1(\Omega,\rho)} \leq \|g_F\|_{L^p(\Omega,\rho)} 
\|\Phi(y)\|_{L^q(\Omega,\rho)},$$ using the fact that $\cX = L^q(\Omega,\rho)$.
Therefore, by the essential boundedness of $\Phi$ w.r.t. $\mu$,
the integrals in (\ref{eq:fub}) are at most
$$\|g_F\|_{L^p(\Omega,\rho)} \|\Phi\|_{\cL^{\infty}(Y,\mu;X} 
\|w\|_{L^1(Y,\mu)} < \infty.$$
Hence, $f$ is the Bochner integral of $w \Phi$ w.r.t. $\mu$. 

(iv) We again use Lemma \ref{pr:mvt}.  Let $Y_0$ be a measurable
subset of $Y$ with $\mu(Y_0) = 0$ and for $Y' = Y \setminus Y_0$,
$\Phi(Y') = G$; see the remark following the definition of 
essential supremum.  But restricting $\sgn \circ w$ and $\Phi$ to $Y'$, 
one has $$f = \cB-\int_{Y'} \sgn(w(y)) \Phi(y) d\mu_w(y);$$
hence, $f \in \mu_w(Y) \cl_\cX \conv (\pm G)$.  Thus,
$\|w\|_{L^1(Y,\mu)} = \mu_w(Y) \geq \|f\|_{G,\cX}$.
\eop

\section{An example involving the Bessel potential}

Here we review an example related to the Bessel functions
which was considered in \cite{kks09} for $q=2$. 
In the following section, this Bessel-potential example is used to find 
an inequality related to the Gamma function. 
 
Let $\cF$ denote the Fourier transform, given for $f \in L^1(\bR^d,\lambda)$ 
and $s \in \bR^d$ by
$$\hat{f}(s) = \cF(f)(s) = (2 \pi)^{-d/2} \int_{\bR^d} f(x) \exp(-i s \cdot x) \; dx, $$
where $\lambda$ is Lebesgue measure and $dx$ means $d\lambda(x)$.  For $r>0$, let
$${\hat {\beta_r}}(s)= (1+\|s\|^2)^{-r/2}\,.$$  
Since the Fourier transform is an isometry of $\cL^2$ onto itself (Parseval's
identity), and ${\hat {\beta_r}}$ is in $\cL^2(\bR^d)$ for $r > d/2$ (which
we now assume), 
there is a unique function $\beta_r$, called the {\em Bessel potential} 
of order $r$, having
${\hat {\beta_r}}$ as its Fourier transform.  See, e.g., \cite[p. 252]{ada75}.
If $1 \leq q <
\infty$ and $r > d/q$, then ${\hat {\beta_r}} \in \cL^q(\bR^d)$
and
\begin{equation}
\hspace{1 in} \|\hat{\beta_r}\|_{\cL^q}  = \pi^{d/2q} \left
(\frac{\Gamma(qr/2 - d/2)}{\Gamma(qr/2)} \right )^{1/q}.
\label{eq:lbeta}
\end{equation} \label{lm:Lbes}

\n
Indeed, by radial symmetry, 
$ (\|{\hat {\beta_r}}\|_{\cL^q})^q =
\int_{\bR^d} (1 + \| x \|^2)^{-qr/2} dx =\omega_d I$, where
 $I = \int_0^{\infty}
(1+\rho^2)^{-qr/2} \rho^{d-1}d\rho$ and
$\omega_d := 2 \pi^{d/2} /\Gamma(d/2)$ is the area of the unit
sphere in $\bR^d$ \cite[p. 303]{co60}. Substituting $\sigma
=\rho^2$ and $d\rho = (1/2) \sigma^{-1/2} d\sigma$, and using
\cite[ p. 60]{ca77}, we find that
$$I = (1/2) \int_0^{\infty}
\frac{\sigma^{d/2 -1}}{(1+ \sigma)^{qr/2}} d\sigma =
\frac{\Gamma(d/2) \Gamma(qr/2 - d/2)}{2\Gamma(qr/2)},$$
establishing (\ref{eq:lbeta}).

For $b>0$, let $\gamma_b: \bR^d \to \bR$ denote the scaled Gaussian
$\gamma_b(x) =  e^{-b\| x\|^2}$.
A simple calculation shows that the $L^q$-norm of $\gamma_b$:
\begin{equation}
\hspace{1.4 in} \|\gamma_b\|_{\cL^q} = (\pi/qb)^{d/2q}. 
\label{eq:normg} \end{equation}
Indeed, using $\int_{-\infty}^\infty \exp(-t^2) dt = \pi^{1/2}$, we obtain:
$$\|\gamma_b\|_{\cL^q}^q = \int_{\bR^d} \exp(-b \|x\|^2)^q dx 
= \left (\int_{\bR} \exp(-qb\,t^2) dt \right )^d = (\pi/qb)^{d/2}.$$
\\
\\
\n We now express the Bessel potential as an integral combination of
Gaussians.  The Gaussians are normalized in $L^q$ and the corresponding
weight function $w$ is explicitly given.  The
integral formula is similar to one in Stein \cite{st70}. By our main
theorem, this is an example of
(\ref{eq:th}) and can be interpreted either as a pointwise integral
or as a Bochner integral.

\begin{proposition}
For $d$ a positive integer, $q \in [1, \infty)$, $r > d/q$, and $s
\in \bR^d$
\begin{displaymath}
\hspace{1 in} \hat{\beta}_{r}(s) = \int_0^{\infty}w_r(t) \gamma_t^o(s) \, dt \,,
\end{displaymath}
where $$\gamma_t^o(s) = \gamma_t(s)/\|\gamma_t\|_{\cL^q}$$ 
and 
$$w_r(t) = (\pi/qt)^{d/2q} \, t^{r/2 - 1}\,
e^{-t}/\Gamma(r/2).$$
\label{pr:calpha}
\end{proposition}

\bop Let
$$I = \int_0^{\infty}t^{r/2-1}\,e^{-t} \, e^{-t\|s\|^2}\,dt.$$
Putting $u = t(1 +\|s\|^2)$ and $dt = du (1
+\|s\|^2)^{-1}$, we obtain 
$$I = (1 + \|s\|^2)^{-r/2}
\int_0^{\infty}u ^{r/2-1}\,e^{-u}\,du = \hat{\beta}_{r}(s)
\Gamma(r/2).$$

\n Using the norm of the Gaussian
(\ref{eq:normg}), we arrive at
$$\hat{\beta}_{r}(s) = I / \Gamma(r/2) = \left (\int_0^{\infty} (\pi/qt)^{d/2q}\,
t^{r/2-1}\,e^{-t}\,  \gamma^o_t(s)dt \right ) {\Large /\; }
\Gamma(r/2),$$ 
which is the result desired. \eop

Now we apply Theorem \ref{th:main} with $Y = (0, \infty)$
and $\phi(s,t) = \gamma_t^o(s) = \gamma_t(s)/\|\gamma_t\|_{L^q(\bR^d)}$
to bound the variational norm of $\hat{\beta_r}$ by the $L^1$-norm of
the weight function.

\begin{proposition} For $d$ a positive integer, $q \in [1, \infty)$,
and $r > d/q$,
$$\|\hat{\beta_r}\|_{G, \cX} \leq (\pi/q)^{d/2q}
\frac{\Gamma(r/2 - d/2q)}{\Gamma(r/2)},$$ where $G =
\{\gamma_t^o : 0 < t < \infty \}$ and $\cX = \cL^q(\bR^d)$.
\label{co:brvar}
\end{proposition}

\bop By (\ref{ineq:varL1}) and Proposition \ref{pr:calpha}, we
have
$$\|\hat{\beta_r}\|_{G, \cX} \leq \|w_r\|_{\cL^1(Y)}
= k \int_0^{\infty} e^{-t} t^{r/2 + d/2q -1} dt,$$ where $k =
(\pi/q)^{d/2q} /\Gamma(r/2)$, and by definition, the integral is $\Gamma(r/2 -
d/2q)$. \eop

\section{Application: A Gamma function inequality}

The inequalities among the variational norm $\|\cdot\|_{G,\cX}$,
the Banach space norm $\|\cdot\|_\cX$, and the $L^1$-norm of the
weight function, established in the Main Theorem,
allow us to derive other inequalities.  The Bessel
potential $\beta_r$ of order $r$ considered above provides an example.

Let $d$ be a positive integer, $q \in [1, \infty)$, and $r > d/q$.  
By Proposition \ref{co:brvar} and (\ref{eq:lbeta}) of the last section, and 
by (\ref{ineq:thm}) of the Main Theorem,
we have
\begin{equation} \hspace{.5 in} 
\pi^{d/2q}\left (\frac{\Gamma(qr/2 - d/2)}{\Gamma(qr/2)}
\right )^{1/q} \leq (\pi/q)^{d/2q} \;\frac{\Gamma(r/2 -
d/2q)}{\Gamma(r/2)}.
\label{eq:gfineq}
\end{equation}
Hence, with $a = r/2 - d/2q\,$ and $s = r/2$, this becomes
\begin{equation}
\hspace{1 in} q^{d/2q} \;\left(\frac{\Gamma(qa)}{\Gamma(qs)}\right)^{1/q} \leq
\frac{\Gamma(a)}{\Gamma(s)}.
\label{eq:g2id}
\end{equation}

In fact, (\ref{eq:g2id}) holds if $s,a,d,q$ satisfy (i) $s > a > 0$ and
(ii) $s - a = d/2q$ for some $d \in Z^+$ and $q \in [1,\infty)$.  As $a>0$,
$r > d/q$. If
$T = \{t > 0 : t = d/2q \;$ for some $\; d \in Z^+, q \in [1,\infty)\}$,
then $T = (0,\frac 1 2] \cup (0,1] \cup (0,\frac 3 2] \cup \ldots = (0,\infty)$,
so there always exist $d$, $q$ satisfying (ii);
the smallest such $d$ is $ \lceil 2 (s - a)\rceil$. 

The inequality (\ref{eq:g2id}) suggests that the Main Theorem can be
used to establish other inequalities of interest among classical functions. 
We now give a direct argument for the inequality.  
Its independent proof confirms our 
function-theoretic methods and provides additional generalization.  

We begin by noting that in (\ref{eq:g2id}) it suffices to take
$d = 2q(s-a)$.  If the inequality is true in that case, it is true
for all real numbers $d \leq 2q(s-a)$.  Thus, we wish to establish that

\begin{eqnarray*}  s \longmapsto \frac{\Gamma(qs)}{\Gamma(s)^q
q^{sq}}\end{eqnarray*}
\noindent is a strictly increasing function of $s$ for $q > 1$ and
$s > 0$. (For $q = 1$ this function is constant.)

Equivalently, we show that
\begin{eqnarray*} H_q(s) := \log{\Gamma(qs)} - q\log{\Gamma(s)}
- sq \log{q} \end{eqnarray*}
\noindent is a strictly increasing function of $s$ for $q > 1$ and
$s > 0$.

Differentiating with respect to $s$, we obtain:
\begin{eqnarray*} \frac{dH_q(s)}{ds} & = & q
\frac{\Gamma'(qs)}{\Gamma(qs)} - q\frac{\Gamma'(s)}{\Gamma(s)} - q
\log{q} \\
& = & q (\psi(qs) - \psi(s) - \log{q}) \\
& =: & q A_s(q)
 \end{eqnarray*}
\noindent where $\psi$ is the digamma function. It suffices to
 establish that $A_s(q) > 0$ for $q > 1$, $s > 0$. Note that $A_s(1)
 = 0$. Now consider
 \begin{eqnarray*} \frac{dA_s(q)}{dq} = s \psi'(qs) - \frac{1}{q}.
 \end{eqnarray*}
\noindent This derivative is positive if and only if $\psi'(qs) >
 \frac{1}{qs}$ for $q > 1$, $s >0 $.

 It remains to show that $\psi'(x) > \frac{1}{x}$ for $x > 0$. Using
 the power series for $\psi'$ \cite[6.4.10]{as}, we have
 for $x > 0$,
 \begin{eqnarray*} \psi'(x)& = & \sum_{n = 0}^{\infty} \frac{1}{(x +
 n)^2} = \frac{1}{x^2} + \frac{1}{(x + 1)^2} + \frac{1}{(x + 2)^2} +
 \ldots  \\ & > &  \frac{1}{x(x + 1)} + \frac{1}{(x + 1)( x + 2)} +
 \frac{1}{(x + 2)(x + 3)} + \dots\\
& = & \frac{1}{x} - \frac{1}{x + 1} + \frac{1}{x + 1} - \frac{1}{x +
2} +
 \ldots \;= \;\frac{1}{x}.
 \end{eqnarray*}

\section{Tensor-product interpretation}

The basic paradigm of feedforward neural nets is to select a single type of computational 
unit and then build a network based on this single type through a choice of 
controlling internal and external parameters so that the resulting network function 
approximates the target function; see \cite{kks-h}.
However, a single type of hidden unit may not be as effective as one 
based on a {\it plurality} of hidden-unit types.  Here we explore a
tensor-product interpretation which may facilitate such a change in
perspective.

Long ago Hille and Phillips \cite[p. 86]{hp57} observed that 
the Banach space of Bochner integrable functions
from a measure space $(Y,\mu)$ into a Banach space $\cX$ has a fundamental set
consisting of two-valued functions, achieving a single non-zero value on
a measurable set of finite measure.  Indeed, every Bochner integrable function
is a limit of simple functions, and each simple function (with a finite
set of values achieved on disjoint parts $P_j$ of the partition) can be written as
a sum of characteristic functions, weighted by members of the Banach space.
If $s$ is such a simple function, then $$s = \sum_{i=1}^n \chi_j g_j,$$  
where the $\chi_j$ are the characteristic functions of the $P_j$ and the $g_j$
are in $\cX$.  
(If, for example, $Y$ is embedded in a finite-dimensional 
Euclidean space, the partition could consist of generalized rectangles.)

Hence, if $f = \cB-\int_Y h(y) d\mu(y)$ is the Bochner integral of $h$ with respect to some measure $\mu$, 
then $f$ can be approximated as closely as desired by elements in $\cX$ of 
the form
$$\sum_{i=1}^n \mu(P_i) g_i,$$
where $Y = \bigcup_{i=1}^n P_i$ is a $\mu$-measurable partition of $Y$.

Note that given a $\sigma$-finite measure space $(Y,\mu)$ and a separable
Banach space $\cX$, every element $f$ in $\cX$ is 
(trivially) the Bochner integral of any 
integrand $w \cdot \kappa(f)$, where $w$ is a nonnegative function on $Y$ with
$\|w\|_{L^1(Y,\mu)} = 1$ (see part (iii) of Theorem \ref{th:main}) 
and $\kappa(f)$ denotes the
constant function on $Y$ with value $f$.  
In effect, $f$ is in $\cX_G$ when $G = \{f\}$.
When $\Phi$ is chosen first (or more precisely $\phi$ as in our Main
Theorem), then $f$ may or may not be in $\cX_G$.  According to the
Main Theorem, $f$ is in $\cX_G$ when it is given by an integral formula
involving $\Phi$ and some $L^1$ weight function.  In this case, 
$G = \Phi(Y) \cap B$ where $B$ is the ball in $\cX$ of radius 
$\|\Phi\|_{L^\infty(Y,\mu;\cX)}$.

In general, the elements $\Phi(y),\; y \in Y$ of the Banach space
involved in some particular approximation for $f$ will be distinct functions
of some general type obtained by varying the parameter $y$.  
For instance, 
kernels, radial basis functions, perceptrons, or various other classes
of computational units can be used, and when these computational-unit-classes determine fundamental sets, by Proposition \ref{pr:varban}, it is
possible to obtain arbitrarily good approximations.  However,
Theorem \ref{th:tensor} below suggests that having a finite set of 
distinct types $\Phi_i: Y \to \cX$ may allow a smaller ``cost'' for
approximation, if we regard 
$$\sum_{i=1}^n\|w_i\|_1 \|\Phi_i\|_\infty$$ 
as the cost of the approximation
$$f = \cB-\int_Y \left (\sum_{i=1}^n w_i \Phi_i \right )(y) d\mu(y).$$
   
We give a brief sketch of the ideas, following Light and Cheney \cite{lc85}.

Let $\cX$ and $\cZ$ be Banach spaces.
Let $\cX \otimes \cZ$ denote the linear space of equivalence classes of
formal expressions 
$$
\sum_{i=1}^n f_i \otimes h_i, \; f_i \in \cX,\; h_i \in \cZ, \; n \in \bN,\;\;
\;\sum_{i=1}^m f'_i \otimes h'_i, \; f'_i \in \cX,\; h'_i \in \cZ, \; m \in \bN,$$
where these expressions are equivalent if for every $F \in \cX^*$
$$\sum_{i=1}^n F(f_i)h_i = \sum_{i=1}^m F(f'_i)h'_i,$$
that is, if the associated operators from $\cX^* \to \cZ$ are identical,
where $\cX^*$ is the algebraic dual of $\cX$.
The resulting linear space $\cX \otimes \cZ$ is called the {\it algebraic
tensor product} of $\cX$ and $\cZ$.
We can extend $\cX \otimes \cZ$ to a Banach space by completing it with
respect to a suitable norm.   
Consider the norm defined for $t \in \cX \otimes \cZ$,
\begin{equation}
\gamma(t) = \inf \left \{ \sum_{i=1}^n \|f_i\|_{\cX} \|h_i\|_\cZ \;: 
\; t = \sum_{i=1}^n f_i \otimes h_i \right \}.
\label{eq:gamma}
\end{equation}
and complete the algebraic tensor product with respect
to this norm; the result is denoted $\cX \otimes_\gamma \cZ$. 

In \cite[Thm. 1.15, p. 11]{lc85}, Light and Cheney showed
that for any measure space $(Y,\mu)$ and any Banach space
$\cX$ the linear map 
$$\Lambda_\cX: L^1(Y,\mu) \otimes \cX \to \cL^1(Y,\mu;\cX)$$ given by 
$$\sum_{i=1}^r w_i \otimes g_i \mapsto \sum_{i=1}^r w_i g_i.$$ 
is well-defined and extends to a map 
$$\Lambda_\cX^\gamma: L^1(Y,\mu) \otimes_\gamma \cX \to L^1(Y,\mu;\cX),$$
which is an isometric isomorphism
of the completed tensor product onto the
space $L^1(Y,\mu;X)$ of Bochner-integrable functions. 

The following theorem extends the function $\Lambda_\cX$ via the natural
embedding $\kappa_\cX$ of $\cX$ into the space of essentially bounded
$\cX$-valued functions defined in section 4. 

\begin{theorem}
Let $\cX$ be a separable Banach space and let $(Y,\mu)$ be a $\sigma$-finite
measure space. Then there exists a continuous linear surjection 
$$e = \Lambda^{\infty,\gamma}_\cX: L^1(Y,\mu) \otimes_\gamma L^\infty(Y,\mu;\cX) 
\to L^1(Y,\mu;\cX).$$
Furthermore, $e$ makes the following diagram commutative:

\begin{equation}
\begin{array}[c]{ccc}
L^1(Y,\mu) \otimes_\gamma \cX \;&\stackrel{a}{\longrightarrow}\;&
L^1(Y,\mu;\cX)\\
\downarrow\scriptstyle{b}&\nearrow\scriptstyle{e}&\downarrow\scriptstyle{c}\\
L^1(Y,\mu) \otimes_\gamma L^\infty(Y,\mu;\cX) \;
&\stackrel{d}{\longrightarrow}&
L^1(Y,\mu;L^\infty(Y,\mu;\cX))
\end{array}
\label{eq:diag}
\end{equation}
where the two horizontal arrows $a$ and $d$ are the isometric isomorphisms
$\Lambda_\cX^\gamma$ and $\Lambda_{L^\infty(Y,\mu;\cX)}^\gamma$; 
the left-hand vertical arrow $b$ is induced by $1 \otimes \kappa_\cX$, while
the right-hand vertical arrow $c$ is induced by post-composition with 
$\kappa_\cX$, i.e., for any $h$ in $L^1(Y,\mu;\cX)$,
$$c(h) = \kappa_\cX \circ h: Y \to L^{\infty}(Y,\mu; \cX).$$
\label{th:tensor}
\end{theorem}

\bop
The map $$e':\sum_{i=1}^n w_i \otimes \Phi_i  \mapsto \sum_{i=1}^n w_i \Phi_i$$
defines a linear function $L^1(Y,\mu) \otimes L^\infty(Y,\mu;\cX) 
\to L^1(Y,\mu;\cX)$; indeed, it 
takes values in the Bochner integrable functions as by our Main Theorem
each summand is in the class.  

To see that $e'$ extends to $e$ on the
$\gamma$-completion, 
\begin{eqnarray*}
\|e(t)\|_{L^1(Y,\mu;\cX)} & = & \|\cB-\int_Y e(t) d\mu(y) \|_\cX 
\leq \int_Y \|e(t)\|_\cX d\mu(y) = 
\end{eqnarray*}

\begin{eqnarray*}
\int_Y \|\sum_i w_i(y) \Phi_i(y)\|_\cX d\mu(y)  & \leq & 
\int_Y \sum_i |w_i(y)| \|\Phi_i(y)\|_\cX d\mu(y)
\end{eqnarray*} 
\begin{eqnarray*}
\leq \sum_i \|w_i\|_1 \|\Phi_i\|_\infty
\end{eqnarray*}
Hence, $\|e(t)\|_{L^1(Y,\mu;\cX)} \leq \gamma(t)$, so the map $e$ is continuous.
\eop

\section{An example involving bounded variation on an interval}

The following example, more elaborate than the one following 
Proposition \ref{pr:varban},
 is treated in part by Barron \cite{ba93} 
and K\accent23urkov\'a \cite{vk02}. 

Let $\cX$ be the set of equivalence classes of 
(essentially) bounded Lebesgue-measurable 
functions on $[a,b]$,
$a,b < \infty$, i.e., $\cX = L^\infty([a,b])$, with norm
$\|f\|_\cX := \inf \{M : |f(x)| \leq M \;\mbox{for almost every}\; x \in [a,b] \}$.  Let $G$ be the set of equivalence classes of all
characteristic functions of closed intervals
of the forms $[a,b]$, or $[a,c]$ or $[c,b]$ with $a < c < b$.  These functions
are the restrictions of
characteristic functions of closed half-lines to $[a,b]$.  
The equivalence relation is $f \sim g$ if and only if $f(x) = g(x)$
for almost every $x$ in $[a,b]$ (with respect to Lebesgue measure). 

Let $BV([a,b])$ be the set of all equivalence classes of functions 
on $[a,b]$ with bounded variation; that is, each equivalence class contains
a function $f$ such that the {\it total variation} $V(f,[a,b])$ is finite,
where total variation is the largest possible total movement of a discrete
point which makes a finite number of stops as x varies from $a$ to $b$,
maximized over all possible ways to choose a finite list of intermediate
points, that is,
\begin{eqnarray*}
\hspace{-.85 cm} 
V(f,[a,b]) &:=&
\sup \{\sum_{i=1}^{n-1} |f(x_{i+1}) - f(x_i)| : 
n \geq 1, a \leq x_1 < x_2 < \cdots < x_n \leq b \}.
\end{eqnarray*}
 
In fact, each equivalence class $[f]$ contains exactly one function $f^*$ of
bounded variation that satisfies the {\it continuity conditions}:\\
\\

(i) $f^*$ is right-continuous at $c$ for $ c \in [a,b)$, and\\
\\

(ii) $f^*$ is left-continuous at $b$.\\
\\

\n Moreover, $V(f^*,[a,b]) \leq V(f,[a,b])$ for
all $f \sim f^*$. 

To see this, 
recall that every function $f$ of bounded variation is the difference 
of two nondecreasing functions $f = f_1 - f_2$, and $f_1, f_2$ are necessarily right-continuous except at a countable set. We can take 
$f_1(x) := V(f,[a,x]) + K$, where $K$ is an arbitrary constant, 
and $f_2(x) := V(f,[a,x]) + K - f(x)$ for $x \in [a,b]$. Now
redefine both $f_1$ and $f_2$ 
at countable sets to form $f_1^*$ and $f_2^*$ which satisfy the continuity
conditions and are still nondecreasing on $[a,b]$. 
Then $f^* := f_1^* - f_2^*$ also satisfies the continuity conditions. 
It is easily shown that $V(f^*,[a,b]) \leq V(f,[a,b]).$ 
Since any equivalence class in $\cX$ can contain at most one function satisfying
(i) and (ii) above, it follows that $f^*$ is unique and that $V(f^*,[a,b])$
minimizes the total variation for all functions in the equivalence class.
Recall that $\chi_{[a,b]} = \chi([a,b])$ denotes the characteristic function
of the interval $[a,b]$, etc.

\begin{proposition}
Let $\cX = L^\infty([a,b])$ and let $G$ be the subset 
of characteristic functions\\ $\chi([a,b]), \chi([a,c]), \chi([c,b]), a < c < b$
(up to sets of Lebesgue-measure zero).\\
Then $\cX_G = BV([a,b])$, and 
$$\|[f]\|_\cX \leq \|[f]\|_{G,\cX} \leq 2 V(f^*,[a,b]) + |f^*(a)|,$$
where $f^*$ is the member of $[f]$ satisfying the continuity conditions 
(i) and (ii).
\label{pr:pr}
\end{proposition} 
 
\bop Let 
$C_{G,\cX}$ be the set of equivalence classes of functions of the form 

\begin{equation}
(q-r) \chi_{[a,b]} + \sum_{n=1}^k
 (s_n - t_n) \chi_{[a,c_n]} + \sum_{n=1}^k (u_n - v_n) \chi_{[c_n,b]},
\label{eq:convex}
\end{equation}
where $k$ is a positive integer, $q, r \geq 0$, for $1 \leq n \leq k$,
$ s_n , t_n, u_n, v_n \geq 0,\;$ and 
   $$(q + r) + \sum_{n=1}^k (s_n + t_n + u_n + v_n) = 1.$$
All of the functions so exhibited have bounded variation $\leq 1$
and hence $C_{G,\cX} \subseteq BV([a,b])$.

We will prove that a sequence in $C_{G,\cX}$ converges
in $\cX$-norm to a member of $BV([a,b])$ and this will establish
that $B_{G,\cX}$ is a subset of $BV([a,b])$ and hence that $\cX_G$ is
a subset of $BV([a,b])$.  

Let $\{[f_k]\}$ be a sequence in $C_{G,\cX}$
that is Cauchy in the $\cX$-norm.  Without loss of generality, we pass
to the sequence $\{f^*_k\}$, which is Cauchy in the sup-norm
since $x \mapsto |f^*_k(x) - f^*_j(x)|$ satisfies the continuity
conditions (i) and (ii).  
Thus, $\{f^*_k\}$ converges pointwise-uniformly and in the sup-norm
to a function $f$ on $[a,b]$ also satisfying (i) and (ii) 
with finite sup-norm and whose equivalence class has finite $\cX$-norm.
  
Let $\{x_1, \ldots, x_n \}$ satisfy $a \leq x_1 < x_2 < \cdots < x_n \leq b$.
Then
$$\sum_{i=1}^{n-1} |f^*_k(x_{i+1}) - f^*_k(x_i)| \leq V(f_k^*,[a,b]) 
\leq V(f_k,[a,b]) \leq 1$$
for every $k$
, where {\it par abus de notation} $f_k$ denotes the member of $[f_k]$
satisfying (\ref{eq:convex}).
Letting $k$ tend to infinity and then varying $n$ and $x_1, \ldots, x_n$, 
we obtain $V(f,[a,b]) \leq 1$ and so $[f] \in BV([a,b])$.

It remains to show that everything in $BV([a,b])$ is actually in $\cX_G$.
Let $g$ be a nonnegative 
nondecreasing function on $[a,b]$ satisfying the continuity
conditions (i) and (ii) above.  Given a positive integer $n$, there exists
a positive integer $m \geq 2$ and $a = a_1 < a_2 < \cdots < a_m = b$ such that
$g(a_{i+1}^-) - g(a_i) \leq 1/n$ for $i = 1, \ldots, m-1$. Indeed, 
for $2 \leq i \leq m-1$, let 
$a_i := \min \{x | g(a) + \frac{(i - 1)}{n} \leq g(x) \}$.  (Moreover,
it follows that the
set of $a_i$'s include all points of left-discontinuity of $g$ such that 
the jump $g(a_i) - g(a_i^-)$ is greater than $1/n$.)  
Let $g_n:[a,b] \to \bR$ be defined as follows:
$$g_n := g(a_1) \chi_{[a_1,a_2)} + g(a_2) \chi_{[a_2,a_3)} + \cdots 
+ g(a_{m-1}) \chi_{[a_{m-1},a_m]}$$ 
$$ = g(a_1) (\chi_{[a_1,b]} - \chi_{[a_2,b]})
+ g(a_2) (\chi_{[a_2,b]} - \chi_{[a_3,b]}) + \cdots 
+ g(a_{m-1}) (\chi_{[a_{m-1},b]})$$ 
$$ = g(a_1) \chi_{[a_1,b]} + (g(a_2) - g(a_1))\chi_{[a_2,b]} + \cdots 
+ (g(a_{m-1}) - g(a_{m-2}))\chi_{[a_{m-1},b]}.$$
Then $[g_n]$ belongs to $g(a_{m-1}) C_{G,\cX}$, and
{\it a fortiori} to $g(b) C_{G,\cX}$ as well as of $\cX_G$, and
$\|[g_n]\|_{G,\cX} \leq g(b)$.
Moreover, $\|[g_n] - [g]\|_\cX \leq 1/n$.  Therefore, since
$B_{G,\cX} = \cl_\cX(C_{G,\cX})$, $[g]$ is in $g(b) B_{G,\cX}$
and accordingly $[g]$ is in $\cX_G$ and $$\|\,[g]\,\|_{G,\cX} \leq g(b).$$

Let $[f]$ be in $BV([a,b])$ and let $f^* = f^*_1 - f^*_2$, 
as defined above, for this purpose we take $K = |f^*(a)|$.
This guarantees that both $f^*_1$ and $f^*_2$ are nonnegative. 
Accordingly, $[f] = [f^*_1] - [f^*_2]$, and is in $ \cX_G$.
Furthermore, 
$\|[f]\|_{G,\cX} \leq \|[f^*_1]\|_{G,\cX} + \|[f^*_2]\|_{G,\cX} \leq f^*_1(b) + f^*_2(b) = V(f^*,[a,b]) + |f^*(a)|
+ V(f^*,[a,b]) + |f^*(a)| - f^*(b) \leq 2V(f^*,[a,b]) + |f^*(a)|$.
The last inequality follows from the fact that 
$V(f^*,[a,b]) + |f^*(a)| - f^*(b) \geq |f^*(b) - f^*(a)| + |f^*(a)| - f^*(b)
\geq 0$. \eop

An argument similar to the above shows that $BV([a,b])$ is a Banach space
under the norm $2 V(f^*,[a,b]) + |f^*(a)|$ (with or without the $2$).
The identity map from $BV([a,b])$ 
(with this norm) to $(\cX_G, \|\cdot\|_{G,\cX}$,
is continuous (by Proposition \ref{pr:pr}) and it is also onto.
Accordingly, by the Open Mapping Theorem (e.g., Yosida \cite[p. 75]{yo65})
the map is open, hence a homeomorphism, so the norms are equivalent.
Thus, in this example, $\cX_G$ is a Banach space under these two equivalent
norms.  

Note however that the $\cX$-norm restricted to $\cX_G$ does not
give a Banach space structure; i.e., $\cX_G$ is not complete in the $\cX$-norm.
Indeed, with $\cX = L^\infty([0,1])$. Let $f_n$ be $1/n$ times the
characteristic function of the disjoint union of $n^2$ closed intervals
contained within the unit interval.  Then $\|[f_n]\|_\cX = 1/n$ 
but $\|[f_n]\|_{G,\cX} \geq C n$, some $C > 0$, since the $\|\cdot\|_{G,\cX}$
is equivalent to the total-variation norm.  While $\{f_n\}$ converges to zero
in one norm, in the other it blows up.  If $\cX_G$ were a Banach space under
$\|\cdot\|_\cX$, it would be another Cauchy sequence, a contradiction.



\section{Pointwise-integrals vs. Bochner integrals} 

\subsection*{Evaluation of Bochner integrals}

A natural conjecture is that the Bochner integral, evaluated pointwise,
is the pointwise integral; that is, if 
$h \in \cL^1(Y,\mu,\cX)$, where $\cX$ is any Banach space of functions
defined on a measure space $\Omega$, then
\begin{equation}
\left (\cB-\int_Y h(y)\, d\mu(y) \right )(x) = \int_Y h(y)(x)\, d\mu(y)
\label{eq:nat}
\end{equation}
for all $x \in \Omega$.  Usually, however, one is dealing with equivalence
classes of functions and thus can expect the equation (\ref{eq:nat}) to
hold only for almost every $x$ in $\Omega$.  Furthermore, to specify
$h(y)(x)$, it is necessary to take a particular function representing 
$h(y) \in \cX$

The Main Theorem implies that (\ref{eq:nat}) holds for $\rho$-a.e. 
$x \in \Omega$
when $\cX = L^q(\Omega, \rho)$, for $1 \leq q < \infty$, is separable 
provided that $h = w \Phi$, where $w: Y \to \bR$ is a
weight function with finite $L^1$-norm and $\Phi: Y \to \cX$ is essentially bounded, where for each $y \in Y$
$\Phi(y)(x) = \phi(x,y)$ for $\rho$-a.e. $x \in \Omega$
and $\phi: \Omega \times Y \to \bR$
is $\rho \times \mu$-measurable.  More generally, we can show the following.

\begin{theorem}
Let $(\Omega,\rho)$, $(Y,\mu)$ be $\sigma$-finite measure spaces, let  
$\cX = L^q(\Omega,\rho)$, $q \in [1,\infty]$, and
let $h \in \cL^1(Y,\mu;\cX)$ so that for each $y$ in $Y$,
$h(y)(x) = H(x,y)$ for $\rho$-a.e. $x$, where $H$ is a
$\rho \times \mu$-measurable real-valued function on $\Omega \times Y$.  Then\\
(i) $y \mapsto H(x,y)$ is integrable for $\rho$-a.e. $x \in \Omega$,\\
(ii) the equivalence class of $x \mapsto \int_Y H(x,y)\, d\mu(y) \mbox{ is  in }\cX$,
and\\
(iii) for $\rho$-a.e. $x \in \Omega$
$$\left (\cB-\int_Y h(y)\, d\mu(y) \right )(x) = \int_Y H(x,y)\, d\mu(y).$$
\label{th:eval}
\end{theorem}

\bop We first consider the case $1 \leq q < \infty$.
Let $g$ be in $\cL^p(\Omega,\rho)$, where $1/p + 1/q = 1$.  Then 
\begin{eqnarray}
\int_Y \int_\Omega |g(x) H(x,y)| d\rho(x) d\mu(y)  
\leq \int_Y \|g\|_p \|h(y)\|_q \,d\mu(y)\\ 
= \|g\|_p \int_Y \|h(y)\|_\cX \,d\mu(y) < \infty.
\label{eq:hs}
\end{eqnarray}
Here we have used Young's inequality and Bochner's theorem.
By Fubini's theorem, (i) follows.  In addition, the map 
$g \mapsto \int_\Omega g(x)\left (\int_Y H(x,y) \,d\mu(y) \right) d\rho(x)$
is a continuous linear functional $F$ on $L^p$ with 
$\|F\|_{\cX^*} \leq \int_Y \|h(y)\|_\cX \,d\mu(y)$.
Since $(L^p)^*$ is $L^q$ for $1 < q < \infty$, then the function 
$x \mapsto \int_Y H(x,y) \,d\mu(y)$ is in $\cX = L^q$ and has norm
$\leq \int_Y \|h(y)\|_\cX \,d\mu(y)$.
The case $q=1$ is covered by taking $g \equiv 1$, a member of $L^\infty$,
and noting that $\|\int_Y H(x,y) \,d\mu(y)\|_{L^1} 
=  \int_\Omega | \int_Y H(x,y) \,d\mu(y) | d\rho(x)
\leq \int_\Omega \int_Y |H(x,y)| \,d\mu(y) d\rho(x) 
= \int_Y \|h(y)\|_\cX \,d\mu(y)$.
Thus, (ii) holds for $1 \leq q < \infty$.

Also by Fubini's theorem and Theorem \ref{pr:blf}, for all $g \in \cX^*$,
$$ \int_\Omega g(x) \left (\cB-\int_Y h(y) \,d\mu(y) \right )(x) d\rho(x) =
\int_Y \left (\int_\Omega g(x)  H(x,y)  d\rho(x) \right ) d\mu(y) $$
$$=
\int_\Omega g(x) \left (\int_Y H(x,y) d\mu(y)\right ) d\rho(x).$$
Hence (iii) holds for all $q < \infty$, including $q = 1$. 

Now consider the case $q = \infty$. For 
$g \in L^1(\Omega,\rho) = \left(L^\infty(\Omega,\rho)\right)^*$, 
the inequality (\ref{eq:hs}) holds, and by \cite[pp. 348--9]{hs65},
(i) and (ii) hold and 
$\|\int_Y H(x,y) \,d\mu(y)\|_\infty \leq \int_Y \|h(y)\|_\infty \,d\mu(y) < \infty$. 
For $g \in \cL^1(\Omega, \rho)$,
$$\int_\Omega g(x) \left (\cB-\int_Y  h(y) \,d\mu(y) \right )(x) d\rho(x) =
\int_Y \left ( \int_\Omega g(x) h(y)(x) d\rho(x) \right ) d\mu(y)$$ 
$$= \int_Y \left ( \int_\Omega g(x) H(x,y) d\rho(x) \right ) d\mu(y)
= \int_\Omega g(x) \left (\int_Y  H(x,y) \,d\mu(y) \right ) d\rho(x).$$
The two functions integrated against $g$ are in $L^\infty(\Omega, \rho)$
and agree, so the functions must be the same $\rho$-a.e. \eop 
 
There are cases where $\cX$ consists of pointwise-defined functions
and (\ref{eq:nat}) can be taken literally.

If $\cX$ is a separable Banach space of pointwise-defined functions from
$\Omega$ to $\bR$ in which the evaluation functionals are bounded
(and so in particular if $\cX$ is a reproducing kernel Hilbert space \cite{ar50}),
then (\ref{eq:nat}) holds for all $x$ (not just $\rho$-a.e.).  Indeed, for each
$x \in \Omega$, the evaluation functional $E_x: f \mapsto f(x)$ is bounded and linear,
so by Theorem \ref{pr:blf}, $E_x$ commutes with the Bochner integral operator.
As non-separable reproducing kernel Hilbert spaces exist \cite[p.26]{da01}, one still needs the hypothesis of separability.  

In a special case involving Bochner integrals with values in
Marcinkiewicz spaces,  Nelson \cite{rn82} showed that
(\ref{eq:nat}) holds. His result involves going from equivalence
classes to functions, and uses a ``measurable selection.''
Reproducing kernel Hilbert spaces were studied by Le Page in \cite{rl72} 
who showed that (\ref{eq:nat}) holds when $\mu$ is a probability measure on $Y$
under a Gaussian distribution assumption on variables in the dual
space.  Another special case of (\ref{eq:nat}) is derived in 
Hille and Phillips \cite[Theorem 3.3.4, p. 66]{hp57}, where
the parameter space is an interval of the real line and the
Banach space is a space of bounded linear transformations (i.e., the
Bochner integrals are operator-valued).

\subsection*{Essential boundedness is needed for the Main Theorem}

The following is an example of a function $h:Y \to \cX$ which is not
Bochner integrable.
Let $Y = (0,1) = \Omega$ with $\rho = \mu =\;$ Lebesgue measure
and $q=1=d$ so $\cX=L^1((0,1))$.  Put $h(y)(x) = y^{-x}$.  Then for all $y \in (0,1)$
$$\|h(y)\|_X = \int_0^1  y^{-x} dx = \frac{1 - \frac 1 y}{\log y}.$$
By l'Hospital's rule 
$$\lim_{y \to 0^+} \|h(y)\|_X = +\infty.$$
Thus, the function $y \mapsto \|h(y)\|_X$ is not essentially bounded 
on $(0,1)$ and Theorem \ref{th:main} does not apply.  Furthermore,
for $y \leq 1/2$, $$\|h(y)\|_\cX \geq \frac {1} {-2y \log y}$$
and $$\int_0^1 \|h(y)\|_\cX \,dy \geq \int_0^{1/2} \frac {1} {-2y \log y} \,dy
= -(1/2) \log{(\log y)}|_0^{1/2} = \infty.$$
Hence, by Theorem \ref{th:bo}, $h$ is not Bochner integrable.
Note however that 
$$f(x) =  \int_Y h(y)(x) d\mu(y) = \int_0^1 y^{-x} dy = \frac {1} {1-x}$$ for every $x \in \Omega$.
Thus $h(y)(x)$ has a pointwise integral $f(x)$ for all $x \in (0,1)$,
but $f$ is not in $\cX = L^1((0,1))$. 

\subsection*{Connection with sup norm}

In \cite{kavkvo07}, we take $\cX$ to be the space of bounded measurable
functions on $\bR^d$, $Y$ equal to the product $S^{d-1} \times \bR$ with
measure $\nu$ which is the (completion of the) product measure determined by
the standard (unnormalized) measure $d(e)$ on the sphere and ordinary Lebesgue measure on $\bR$.  We take $\phi(x,y) := \phi(x,e,b) = \vartheta(e \cdot x + b)$,
so $x \mapsto \vartheta(e \cdot x + b)$ is the characteristic function of the
closed half-space $\{x: e \cdot x + b \geq 0 \}$.

We showed that if a function $f$ on $\bR^d$ 
decays, along with its partials of order $\leq d$, at a sufficient rate,
then there is an integral formula expressing $f(x)$ as an integral combination
of the characteristic functions of closed half-spaces weighted by iterated Laplacians integrated over half-spaces.  The characteristic functions all have
sup-norm of 1 and the weight-function is in $L^1$ of $(Y,\nu)$, where
$Y = S^{d-1} \times \bR$ and $\nu$ is the (completion of the)
product measure determined by
the standard (unnormalized) measure $d(e)$ on the sphere $S^{d-1}$ of unit
vectors in $\bR^d$ and ordinary Lebesgue measure on $\bR$.  

For example, when $d$ is odd, 

$$f(x) = \int_{S^{d-1} \times \bR} w_f(e,b) \vartheta(e \cdot x + b) d\nu(e,b),$$

\n where 

$$w_f(e,b) := a_d \int_{H_{e,b}} D_e^{(d)}f(y) d_H(y),$$

\n with $a_d$ a scalar exponentially decreasing with $d$.
The integral is of the iterated directional derivative over the
hyperplane with normal vector $e$ and offset $b$,
 
$$H_{e,b} := \{y \in \bR^d : e \cdot y + b = 0 \}.$$ 

For $\cX = \cM(\bR^d)$, the space of bounded Lebesgue-measurable functions on
$\bR^d$, which is a Banach space w.r.t. sup-norm, and $G$ the family $H_d$ consisting of the set of all characteristic functions
for closed half-spaces in $\bR^d$, it follows from Theorem \ref{th:main}
that $f \in \cX_G$.

Hence, from the Main Theorem, 
$$f = \cB-\int_{S^{d-1} \times \bR} w_f(e,b) \Theta(e,b)\, d\nu(e,b)$$ 
is a Bochner integral, where $\Theta(e,b) \in \cM(\bR^d)$ is given by

$$\Theta(e,b)(x) := \vartheta(e \cdot x + b).$$

Application of the Main Theorem requires only that $w_f$ be in $L^1$,
but \cite{kavkvo07} gives explicit formulas for $w_f$ (in both
even and odd dimensions) provided that $f$ satisfies the decay conditions
described above and in our paper; see also the other chapter in this book 
referenced earlier.

\section{Some concluding remarks}

Neural networks express a function $f$ 
in terms of a combination of members of a given family $G$ of functions.  
It is reasonable to expect that a function $f$ can be so represented
if $f$ is in $\cX_G$.  The choice of $G$ thus dictates the $f$'s
that can be represented (if we leave aside what combinations are permissible).
Here we have focused on the case $G = \{\Phi(y) : y \in Y \}$.
The form $\Phi(y)$ is usually associated with a specific family such as Gaussians
or Heavisides.  The tensor-product interpretation suggests the
possibility of using multiple families $\{\Phi_j: j \in J \}$ or multiple $G$'s to 
represent a larger class of $f$'s.  Alternatively, one may replace $Y$ by
$Y \times J$ with a suitable extension of the measure.

The Bochner integral approach also permits $\cX$ to be an arbitrary Banach space
(not necessarily an $L^p$-space). For example,
if $\cX$ is a space of bounded linear transformations and $\Phi(Y)$ is
a family of such transformations, we can approximate other members
$f$ of this Banach space $\cX$ in a neural-network-like manner.  Even more
abstractly, we can approximate an {\it evolving} function $f_t$, where
$t$ is time, using weights that evolve over time and/or a family $\Phi_t(y)$
whose members evolve in a prescribed fashion.  Such an approach would require
some axiomatics about permissible evolutions of $f_t$, perhaps similar to methods
used in time-series analysis and stochastic calculus.  See, e.g., \cite{ab82}.

Many of the restrictions we have imposed in earlier sections are not truly
essential.  For example, the separability constraints can be weakened. Moreover,
$\sigma$-finiteness of $Y$ need not be required
since an integrable function $w$ on $Y$ must vanish outside a $\sigma$-finite
subset.  More drastically, the integrable function $w$ can be replaced by a
distribution or a measure.  Indeed,
we believe that both finite combinations and integrals can be
subsumed in generalized combinations derived from Choquet's theorem.
The abstract transformations of the concept of neural network discussed here
provide an ``enrichment'' that may have practical consequences.

\section{Appendix I: Some Banach space background}

The following is a brief account of the machinery of functional analysis
used in this chapter.  See, e.g., \cite{yo65}.
For $G \subseteq \cX$, with $\cX$ any linear space, let 
$$\span_n(G) := \left \{x \in \cX : \exists w_i \in \bR, 
g_i \in G, \; 1 \leq i \leq n, \; \ni \;
x = \sum_{i=1}^n w_i g_i \; \right \}$$ 
denote the set of all $n$-fold linear combinations from $G$. If the
$w_i$ are non-negative with sum $1$, then the combination is called
a {\it convex} combination; $\conv_n(G)$ denotes the set of all
$n$-fold convex combinations from $G$.  Let 
$$\span(G) := \bigcup_{n=1}^\infty \, \span_n(G) \;\;\mbox{and} \;\;
\conv(G) := \bigcup_{n=1}^\infty \, \conv_n(G).$$

A {\it norm} on a linear space $\cX$ is a function
which associates to each element $f$ of $\cX$ a real number $\|f\| \geq 0$
such that\\ 
\\

(1) $\|f\| = 0$ $\iff$  $f = 0$;\\

(2) $\|r f\| = |r| \|f\|$ for all $r \in \bR;\; \mbox{and}$\\

(3) the triangle inequality holds: 
$\|f + g\| \leq \|f\| + \|g\|,\; \forall f,g \in \cX$.\\

\n A metric $d(x,y) := \|x - y\|$ is defined by the norm, and both
addition and scalar multiplication become continuous functions
with respect to the topology induced by the norm-metric.
A metric space is {\it complete} if
every sequence in the space that satisfies the Cauchy criterion is convergent.
In particular, if a normed linear space is complete in the metric
induced by its norm, then it is called a {\it Banach} space. 

Let $(Y,\mu)$ be a measure space; it is called $\sigma$-finite provided
that there exists a countable family $Y_1, Y_2, \ldots$ of subsets of $Y$ 
pairwise-disjoint
and measurable with finite $\mu$-measure such that $Y = \bigcup_i Y_i$.
The condition of $\sigma$-finiteness is required for Fubini's theorem.
A set $N$ is called a $\mu$-null set if it is measurable with $\mu(N) = 0$.
A function from a measure space to another measure space
is called {\it measurable} if the pre-image of each measurable
subset is measurable.  When the range space is merely a topological
space, then functions are measurable if the pre-image of each open
set is measurable.

Let $(\Omega,\rho)$ be a measure space.
If $q \in [1,\infty)$, we write $L^q(\Omega, \rho)$ for the
Banach space consisting of all equivalence classes of the set
$\cL^q(\Omega, \rho)$ of all $\rho$-measurable functions from $\Omega$
to $\bR$ with absolutely integrable $q$-th powers,
where $f$ and $g$ are equivalent if they agree $\rho$-almost everywhere
($\rho$-a.e.) - that is, if the set of points where $f$ and $g$ disagree
has $\rho$-measure zero, and 
$\|f\|_{L^q(\Omega, \rho)} := (\int_{\Omega} |f(x)|^q d\rho(x))^{1/q}$,
or $\|f\|_q$ for short. 

\section{Appendix II: Some key theorems}

We include, for the reader's convenience, the statements of some crucial
theorems cited in the text.  

The following consequence of the Hahn-Banach
Theorem, due to Mazur, is given by Yosida \cite[Theorem 3', p. 109]{yo65}.
The hypotheses on $\cX$ are satisfied by any Banach space, but the 
theorem holds much more generally.  See \cite{yo65} for examples where
$\cX$ is not a Banach space.

\begin{theorem}
Let $X$ be a real locally convex linear topological space, $M$ a closed convex
subset, and $x_0 \in X \setminus M$.  
Then $\exists$ continuous linear functional 
$$F: X \to \bR \; \ni \; F(x_0) > 1,\; F(x) \leq 1 \;\; \forall x \in M.$$
\label{th:mazur}
\end{theorem}

\smallskip

Fubini's Theorem relates iterated integrals to product integrals. Let 
$Y, Z$ be sets and $\cM$ be a $\sigma$-algebra of subsets of $Y$ and $\cN$
a $\sigma$-algebra of subsets of $Z$.  If $M \in \cM$ and $N \in \cN$, then
$M \times N \subseteq Y \times Z$ is called a {\it measurable rectangle}.
We denote the smallest $\sigma$-algebra on $Y \times Z$ which contains all
the measurable rectangles by $\cM \times \cN$.  Now let $(Y,\cM,\mu)$ and
$(Z,\cN,\nu)$ be $\sigma$-finite measure spaces, and for $E \in \cM \times \cN$,
define 
$$(\mu \times \nu)(E) := \int_Y \nu(E_y) d\mu(y) = \int_Z \mu(E^z) d\nu(z),$$
where $E_y := \{z \in Z : (y,z) \in E \}$ and 
$E^z := \{y \in Y : (y,z) \in E \}$.  Also, $\mu \times \nu$
is a $\sigma$-finite measure on $Y \times Z$ with $\cM \times \cN$ as the
family of measurable sets.  For the following, see 
Hewitt and Stromberg \cite[p. 386]{hs65}.

\begin{theorem}
Let $(Y,\cM,\mu)$ and $(Z,\cN,\nu)$ be $\sigma$-finite measure spaces.
Let $f$ be a complex-valued $\cM \times \cN$-measurable function on $Y \times Z$,
and suppose that at least one of the following three absolute integrals is finite:
$\int_{Y \times Z} |f(y,z)| d(\mu \times \nu)(y,z)$,
$\int_Z \int_Y |f(y,z)| d\mu(y) d\nu(z)$,
$\int_Y \int_Z |f(y,z)| d\nu(z) d\mu(y)$.
Then the following statements hold:\\
(i) $y \mapsto f(y,z)$ is in $\cL^1(Y,\cM,\mu)$ for $\nu$-a.e. $z \in Z$;\\
(ii) $z \mapsto f(y,z)$ is in $\cL^1(Z,\cN,\nu)$ for $\mu$-a.e. $y \in Y$;\\
(iii) $z \mapsto \int_Y f(y,z) d\mu(y)$ is in $\cL^1(Z,\cN,\nu)$;\\
(iv) $y \mapsto \int_Z f(y,z) d\nu(z)$ is in $\cL^1(Y,\cM,\mu)$;\\
(v) all three of the following integrals are equal:
$$\int_{Y \times Z} f(y,z) d(\mu \times \nu)(y,z) = $$
$$\int_Z \int_Y f(y,z) d\mu(y) d\nu(z) = $$
$$\int_Y \int_Z f(y,z) d\nu(z) d\mu(y).$$
\label{th:fubini}
\end{theorem}

\smallskip

A function $G:I \to \bR$, $I$ any subinterval of $\bR$, 
is called {\it convex} if 
$$\forall x_1, x_2 \in I, 0 \leq t \leq 1, \;\;
G(t x_1 + (1-t) x_2) \leq t G(x_1) + (1-t) G(x_2).$$
The following formulation is from Hewitt and Stromberg \cite[p. 202]{hs65}.

\begin{theorem}[Jensen's inequality]
Let $(Y,\sigma)$ be a probability measure space.  Let $G$ be
a convex function from an interval $I$ into $\bR$ and let 
$f$ be in $\cL^1(Y,\sigma)$ with $f(Y) \subseteq I$ such that
$G \circ f$ is also in $\cL^1(Y,\sigma)$.  Then $\int_Y f(y) d\sigma(y)$
is in $I$ and
$$G \left (\int_Y f(y) d\sigma(y) \right ) \leq \int_Y (G \circ f)(y) 
d\sigma(y).$$
\label{th:jensen}
\end{theorem}

\section*{Acknowledgements}

We thank Victor Bogdan for helpful comments on earlier versions.

\end{document}